\documentclass{amsart}
\usepackage{graphicx}
\usepackage{latexsym}
\usepackage{hyperref}

\usepackage{amsfonts}
\usepackage[all]{xy}

\numberwithin{equation}{section}
\newtheorem{thm}{Theorem}[section]
\newtheorem{cor}[thm]{Corollary}
\newtheorem{lem}[thm]{Lemma}

\theoremstyle{definition}

\numberwithin{equation}{section} \theoremstyle{remark}

\renewcommand{\theequation}{\arabic{section}.\arabic{equation}}

\begin{document}
\renewcommand{\theequation}{\arabic{section}.\arabic{equation}}
\newcommand{\be}{\begin{eqnarray}}
\newcommand{\en}{\end{eqnarray}}
\newcommand{\no}{\nonumber}
\newcommand{\la}{\lambda}
\newcommand{\laa}{\Lambda}
\newcommand{\ino}{\int_\Omega}
\newcommand{\tg}{\triangle}
\newcommand{\Li}{\Lambda_i}
\newcommand{\Lj}{\Lambda_j}
\newcommand{\Lk}{\Lambda_{k+1}}
\newcommand{\li}{\lambda_i}
\newcommand{\lj}{\lambda_j}
\newcommand{\lk}{\lambda_{k+1}}
\newcommand{\p}{\partial}
\newcommand{\n}{\nabla}
\newcommand{\nj}{\nabla_j}
\newcommand{\nk}{\nabla_k}
\newcommand{\mc}{\mcite}
\newcommand{\g}{\gamma}
\newcommand{\e}{\varepsilon}
\newcommand{\s}{\sup}
\newcommand{\ep}{\epsilon}
\newcommand{\de}{\delta}
\newcommand{\D}{\Delta}
\newcommand{\pl}{\parallel}
\newcommand{\ov}{\overline}
\newcommand{\bet}{\beta}
\newcommand{\al}{\alpha}
\newcommand{\fr}{\frac}
\newcommand{\pa}{\partial}
\newcommand{\we}{\wedge}
\newcommand{\om}{\Omega}
\newcommand{\na}{\nabla}
\newcommand{\lan}{\langle}
\newcommand{\ra}{\rangle}
\newcommand{\fa}{\sum_{\alpha=1}^n}
\newcommand{\vs}{\vskip0.3cm}
\newcommand{\su}{\sum_{i=1}^k(\Lambda_{k+1}-\Lambda_i)^2}
\newcommand{\suu}{\sum_{i=1}^k(\Lambda_{k+1}-\Lambda_i)}
\renewcommand{\thefootnote}{}
\newcommand{\ri}{\rightarrow}
\title[] {Pohozaev
identity for the anisotropic $p$-Laplacian and estimates of torsion
function}
 \footnotetext{This work was partially supported by CNPq-BR grant numbers: 307089/2014-2 and 306146/2014-2.}
\author[]
{ Qiaoling Wang and  Changyu Xia}
\address{Qiaoling Wang \\ \newline \indent
Departamento de Matem\'{a}tica, Universidade de Bras\'{\i}lia,
70910-900-Bras\'{\i}lia-DF, Brazil. wang@mat.unb.br}

\address{Changyu Xia \\  \newline \indent   Departamento de Matem\'atica, Universidade de Bras\'{\i}lia, 70910-900 Bras\'{\i}lia-DF, Brazil.  xia@mat.unb.br}
\date{}
\maketitle
\begin{abstract} \noindent In this paper we prove the Pohozaev identity for
the weighted anisotropic $p$-Laplace operator. As an application of
our identity, we deduce the nonexistence of nontrivial solutions of
the Dirichlet problem for the weighted anisotropic $p$-Laplacian in
star-shaped domains of $\mathbb{R}^n$. We also provide an upper
bound estimate for the first Dirichet eigenvalue of the anisotropic
$p$-Laplacian on bounded domains of $\mathbb{R}^n$, some sharp
estimates for the torsion function of compact manifolds with
boundary and a nonexistence result for the solutions of the Laplace
equation on closed Riemannian manifolds.

\end{abstract}


\section{Introduction and the main results}
 Let $\om$ be a bounded smooth domain in $\mathbb{R}^n$ and $g$ a continuous function on $\mathbb{R}$.
In 1965, Pohozaev \cite{Po} considered the following nonlinear
elliptic
problem : \be\label{inteq1}\left\{\begin{array}{l}\label{ineq1} -\D u = g(u) \ \ \ {\rm in\ \ }\om, \\
\ \ \ \ \   u=0  \ \ \ \ \ \ \ {\rm on\ \ } \pa\om,
\end{array}\right.
\en and proved that if $u\in C^2(\om)\cap C^1(\ov{\om})$ is a
solution of (\ref{inteq1}), then \be\label{inteq2} (2-n)\int_{\om}
ug(u)dx +2n \int_{\om} G(u) dx=\int_{\pa\om}|\na u|^2\langle x,
\nu\rangle ds, \en where $G(u)=\int_0^u g(t)dt, \nu=\nu(x)$ is the
outward unit normal vector at the point $x\in \pa\om$.

Based on (\ref{inteq2}), Pohozaev established the following
well-known non-existence result: \vskip0.2cm
 {\bf Theorem A}
(Pohozaev). {\it Let $\om$ be a star-shaped domain with respect to
the origin in $\mathbb{R}^n, n\geq 3$ and $g\in C(\mathbb{R},
\mathbb{R})$ with $ g(u)\geq 0,$ when $u\geq 0$. If
\be\label{inteq3} (2-n)u g(u)+ 2n\int_0^u g(t)dt\leq 0, \ \ {\rm
when} \ u\geq 0,\en then the problem (\ref{inteq1}) has no positive
solution.} \vskip0.2cm

Pohozaev's identity  has also other important applications to the
solutions of differential equations. As an example, let us assume
further that \be g(u)\equiv 1 \ \ {\rm  and}\ \ \left.\fr{\pa
u}{\pa\nu}\right|_{\pa\om}=c=const.\en Then we have from
(\ref{inteq2}) that \be\label{inteq4} (n+2)\int_{\om} u dx&=&c^2
\int_{\pa\om} \langle x, \nu\rangle ds = nc^2 V(\om). \en Also, it
is easy to see in this case that \be \label{inteq5} \D \left(|\na
u|^2 +\fr 2n u\right)=2\left( |\na^2 u|^2 -\fr 1n\right)=2\left(
|\na^2 u|^2 -\fr{(\D u)^2}n\right)\geq 0 \en and \be\label{ineq6}
\left.\left(|\na u|^2 +\fr 2n u\right)\right|_{\pa\om}= c^2, \en
which, implies by the maximum principle that \be \label{inteq7} |\na
u|^2 +\fr 2n u \leq c^2\ \ \ {\rm on\ \ } \om. \en On the other
hand, one obtains from integration by parts and (\ref{inteq4}) that
\be\label{inteq8} \int_{\om}\left(|\na u|^2 +\fr 2n u \right) dx
=\fr{n+2}n \int_{\om} u dx = c^2 V(\om). \en Therefore $|\na u|^2
+\fr 2n u $ is constant in $\om$ and so the equality must hold in
(\ref{inteq5}) which implies that \be \label{inteq9} u_{ij}= -\fr 1n
\delta_{ij} \en Hence, for a suitable choice of origin we know that
$u$ is given by \be u=\fr 1{2n}(\rho_0 -r^2),\en where $\rho_0$ is a
constant and $r$ is the distance function from the origin. Since
$u|_{\pa\om}=0$, we conclude that $\rho_0>0$ and that $\om$ is a
ball of radius $\sqrt{\rho_0}$. Also one can deduce from
(\ref{inteq8}) that $\rho_0=n^2c^2$. The above arguments are
essentially the proof given
 by Weinberger \cite{wein} to the following seminal work of Serrin
\cite{serrin}: \vskip0.2cm {\bf Theorem B} (Serrin). {\it If $u\in
C^2(\ov{\om})$ satisfies the overdetermined problem
\be\left\{\begin{array}{l}\label{serrin1}
\D u = -1 \ \ \ {\rm in\ \ }\om, \\
u|_{\pa \om}=0,
 \left.\fr{\pa
u}{\pa\nu}\right|_{\pa \om}=c,
\end{array}\right.\en where $\om$ is a bounded smooth
domain in $\Bbb R^n$, $\nu$ is the unit outward normal of $\pa \om$,
and $c$ is a constant, then $\om$ is a ball of radius $n|c|$ and
$u=(n^2c^2-r^2)/2n$, where $r$ is the distance from the center of
the ball.} \vskip0.2cm

The appearance of Pohozaev's identity is a milestone in the
developments of differential equations. The generalizations of
Pohozaev's identity have been widely used to prove the non-existence
of nontrivial solutions of nonlinear elliptic equations. Here are
some of the important results in this direction. Esteban-Lions
\cite{EL} and  Berestycki-Lions \cite{BL} considered the following
problem on unbounded domain: \be\left\{\begin{array}{l}
\label{bleq1} -\D u =
g(u) \ \ \ {\rm in\ \ \ }\om,  \\
 u|_{\pa\om}=0,
\end{array}\right. ,\en where $\om =\mathbb{R}^n$ or an unbounded domain of
$\mathbb{R}^n$.  They established the Pohozaev identity for the
above problem and the existence and  nonexistence results which have
brought great developments in this area.
  Pucci and Serrin \cite{PS1} proved the Pohozaev identity satisfied by the general
elliptic equations on bounded domains.  Guedda-Veron \cite{GV}
proved the Pohozaev identity to the solutions of the quasi-linear
elliptic problem \be\left\{\begin{array}{l}\label{eqgv1}
-{\rm div}(|\na u|^{p-2}\na u) = g(x, u) \ \ \ {\rm in\ \ }\om, \\
 u|_{\pa\om}=0,
\end{array}\right.
\en and obtained the non-existence results when $\om$ is a
star-shaped domain. Bartsch-Peng-Zhang \cite{BPZ} and Kou-An
\cite{KA}  considered the more general quasi-linear elliptic
equations with weight on more general domains. Pucci and Serrin
\cite{PS2} studied the Pohozaev identity of polyharmonic operators
and obtained non-existence of nontrivial solutions of the related
equations.

Recently, Ros-Oton and Serra \cite{RS} established the Pohozaev
identity for the fractional elliptic problem:
\be\left\{\begin{array}{l}\label{rseq1}
(-\D )^s u = g(u)  \ \ \ {\rm in\ \ }\om, \\
 u|_{\pa\om}=0 \ \ \ \ \  \ \ \ \ \ \ \ {\rm in\ \ }\mathbb{R}^n\setminus\om
\end{array}\right.
\en in a bounded domain $\om\subset\mathbb{R}^n$, where $s\in (0,
1)$, \be (-\D)^s u(x)= c_{n,s} PV \int_{\mathbb{R}^n} \fr{ u(x)-
u(y)}{|x-y|^{n+2s}} dy \en is the fractional Laplacian and $c_{n,s}$
is a normalization constant given by \be c_{n,s}=\fr{s
2^{2s}\Gamma\left(\fr{n+2s}2\right)}{\pi^{n/2}\Gamma(1-s)}. \en

In this paper, we shall prove the Pohozaev identity for a weighted
anisotropic $p$-Laplace operator.   Let us fix some required
notation before stating our result. Let $F: \mathbb{R}^n\ri [0,
+\infty)$ be a convex function of class
$C^1(\mathbb{R}^n\setminus\{0\})$ which is even and positively
homogeneous of degree $1$, so that \be\label{ineq1} F(tx)=|t|F(x), \
\ \ \forall x\in\mathbb{R}^n, \ \ \ \forall t\in \mathbb{R}. \en
Note that there are positive constants $\alpha$ and $\beta$ such
that $F$ satisfies \be \label{aniseq1} \alpha|\xi |\leq F(\xi)\leq
\beta |\xi|\ \ \ \forall \xi \in \mathbb{R}^n. \en Observing that
$F^p$ is positively homogeneous of degree $p$, we have
\be\label{aniseq3} \langle Z, \na_{\xi}[F^p](Z)\rangle = p F^p(Z), \
\ \forall Z\in \mathbb{R}^n. \en
 For $1<p<\infty$,
the anisotropic $p$-Laplace operator is defined as
\be\label{aniseq2} {\mathcal Q}_p(u)={\rm div}\left(\fr 1p
\na_{\xi}[F^p](\na u)\right)=\sum_{i=1}^n\fr{\pa}{\pa
x_i}\left((F(\na u))^{p-1} F_{\xi_i}(\na u)\right), \en where
$\na_{\xi}$ stands for the gradient operator with respect to the
$\xi$ variables. The first result of the present paper is as
follows.
\begin{thm}\label{th1}
Let $\om$ be a bounded domain with smooth boundary and $g:
\mathbb{R}^n\times\mathbb{R}\ri\mathbb{R}$  a continuous function.
Let  $b$ be a real number, $1<p<\infty$ and $u\in C^2(\om)\cap
C^1(\ov{\om})$ a solution of the problem
\be\label{th1eq1}\left\{\begin{array}{l}
- \fr 1p {\rm div}\left(|x|^{-bp}\na_{\xi}[F^p](\na u)\right)  = g(x, u)  \ \ \ {\rm in\ \ }\om, \\
 u|_{\pa\om}=0.
\end{array}\right.
\en Then we have \be \label{th1eq2} & &\left(1+b-\fr
np\right)\int_{\om} u g(x, u) dx+n\int_{\om} G(x, u) dx
+\int_{\om}\langle x, \na_xG(x, u)\rangle dx\\ \no &=&  \left(1- \fr
1p\right)\int_{\pa\om} F^p(\na u)\langle x, \nu\rangle ds, \en where
$G(x, \rho)=\int_0^\rho g(x,\theta)d\theta$ and $\na_x G$ is the
gradient of $G$ with respect to the first variable $x$.
\end{thm}
As an immediate application of Theorem \ref{th1}, we have
\begin{cor}
Let $\om$ be a bounded star-shaped domain with smooth boundary and
$g: \mathbb{R}\ri\mathbb{R}$  a continuous function. Let $b $ be a
real number, $1<p<\infty$ and suppose that \be
\label{th1eq3}\left(1+b-\fr np \right) \rho g(\rho)+ n\int_0^\rho
g(\sigma)d\sigma < 0, \ \ {\rm when\ } \rho\neq 0. \en Then the
problem \be\label{th1eq4}\left\{\begin{array}{l}
- \fr 1p {\rm div}\left(|x|^{-bp}\na_{\xi}[F^p](\na u)\right)  = g(u)  \ \ \ {\rm in\ \ }\om, \\
 u|_{\pa\om}=0.
\end{array}\right.
\en has no nontrivial (not identically zero) solution.
\end{cor}
Taking \be\label{th1eq5}
 g(x, u)=\lambda |x|^{-\alpha}|u|^{r-2} u
+\mu |x|^{-\beta} |u|^{s-2} u+\eta|x|^{-\gamma} |u|^{t-2} u
 \en
in Theorem\ref{th1}, we have the following
\begin{cor}\label{cor2}
Let $\om$ be a bounded star-shaped domain with smooth boundary and
$b, \alpha, \beta, \gamma, \lambda, \mu, \eta, r, s, t$ be constants
such that $ rst\neq 0, $ and  \be\label{cor2eq1} \la\left(1+b-\fr np
+\fr{n-\alpha}r\right)\leq 0, \ \ \ \mu\left(1+b-\fr
np+\fr{n-\beta}s\right)\leq 0, \\ \no  \hskip0.4cm \eta\left(1+b-\fr
np+\fr{n-\gamma}s\right) < 0. \en Then, the problem
\be\label{th1eq4}\left\{\begin{array}{l}
- \fr 1p {\rm div}\left(|x|^{-bp}\na_{\xi}[F^p](\na u)\right)  = \la |x|^{-\alpha}|u|^{r-2}u+\mu |x|^{-\beta} |u|^{s-2} u  \\
 \hskip4.5cm\ \  + \eta |x|^{-\gamma}|u|^{t-2}u \ \ \ \ \ \ {\rm in\ \
 }\om, \\
u|_{\pa\om}=0,
\end{array}\right.
\en
has no positive solution.
\end{cor}

In the second part of this paper, we study the first Dirichlet
eigenvalue of ${\mathcal Q}_p$ which is given by \be\label{th2eq1}
\la_{p,1}(F,\om)=\inf_{u\in
W_0^{1,p}(\om)\setminus\{0\}}\fr{\int_{\om} F^p(\na u)
dx}{\int_{\om} |u|^p dx}. \en It is known \cite{BFK} that
(\ref{th2eq1}) has a unique positive solution $u_p$ solving the
Euler-Lagrange equation
 \be \label{th2eq2}
\left\{\begin{array}{l}
 {\mathcal Q}_p u_p +\la_{p,1}|u_p|^{p-2}u_p = 0   \ \ {\rm in \ \ } \om, \\
\ \ \ \ \ \ \ \ \ \ \ \ \ \ \ \ \ \ \  \  \ \ \ \ \ u_p = 0  \ \
{\rm on \ \ } \pa \om.
\end{array}\right.
\en The  torsion problem for the anisotropic $p$-Laplace is
 as follows
 \be
\label{th2eq3} \left\{\begin{array}{l}
  -{\mathcal Q}_p v = 1   \ \ {\rm in \ \ } \om, \\
\ \ \ \ \ \ v = 0  \ \ {\rm on \ \ } \pa \om.
\end{array}\right.
\en By classical result there exists a unique solution of
(\ref{th2eq3}), that will be always denoted by $v_{\om}$, which is
positive in $\om $. The anisotropic $p$-torsional rigidity of $\om$
is defined as \be\label{th2eq4} T_{F, p}(\om) = \int_{\om} F^p(\na
v_{\om}) dx =\int_{\om} v_{\om} dx. \en The following variational
characterization for $T_{F, p}(\om)$ holds \be\label{th2eq5}
T_{F,p}(\om)^{p-1} =\max_{\phi\in
W_0^{1,p}(\om)\setminus\{0\}}\fr{\left(\int_{\om}|\phi|dx\right)^p}{\int_{\om}
F^p(\na\phi)dx} \en and the solution $v_{\om}$ of (\ref{th2eq3})
realizes the maximum in (\ref{th2eq5}).

The next result is an estimate involving $\la_{p,1}(F,\om)$ and
$T_{F, p}(\om)$ which is motivated by Theorem 1.1 in \cite{BFNT}.
\begin{thm}\label{th2}Let $p\geq 2$, $F$ as above  and $\om$ a bounded domain in
$\mathbb{R}^n$. We assume further that $F\in
C^{3,\beta}(\mathbb{R}^n\setminus\{0\})$ and \be \label{th201}
[F^p]_{\xi \xi}(\xi) \  is\ positive\ definite \ in\ \
\mathbb{R}^n\setminus\{ 0\}. \en Then, we have \be\label{th2eq7}
 \fr{\la_{p,1}(F,\om)T_{F,p}(\om)^{p-1}}{|\om|^{p-1}}\leq 1 -
\fr{p^{\fr{2p-3}{p-1}}(n\kappa_n^{1/n})^{\fr
p{p-1}}}{n(p-1)(n(p-1)+p)}\cdot \fr{T_{F,p}(\om)}{|\om|^{1+\fr
p{n(p-1)}}}, \en
 where $|\om|$ and $\kappa_n$ stand for the measure
of $\om$ and $K^o$, respectively, being\be\label{th2e1}
K^{o}=\left\{x\in \mathbb{R}^n: \sup_{z\ne 0}\fr{\langle x,
z\rangle}{F(z)}\leq 1\right\}. \en
 \end{thm}
A main tool in the proof of Theorem \ref{th2} is the isoperimetric
inequality (Wulff Theorem) relating the perimeter of a set $E$ with
respect to $F$ and  $|E|$, the measure of $E$. This tool can be also
used to prove the following result which is motivated by \cite{PSS}.

\begin{thm}\label{th0} Let $\om$ be a bounded domain with smooth
boundary in $\mathbb{R}^n$ and $g\in C(\mathbb{R}, \mathbb{R})$ with
$g(\sigma)\geq 0$, when $\sigma\geq 0$. Let $u$ be a smooth positive
solution of the Dirichlet problem for the anisotropic $n$-Laplace
operator: \be\label{th0eq1} \left\{\begin{array}{l}
  -\fr 1n{\rm div}(\na_{\xi}[F^n](\na u)) = g(u)   \ \ {\rm in \ \ } \om, \\
\ \ \ \ \ \ u = 0  \ \ {\rm on \ \ } \pa \om.
\end{array}\right.
\en Then we have \be\label{th0e2} \left(\int_{\om} g(u)
dx\right)^{\fr n{n-1}}\geq \fr{n^{\fr{2n-1}{n-1}}\kappa_n^{\fr
1{n-1}}}{n-1}\int_{\om} G(u) dx,
 \en
where $G(u)=\int_0^u g(s)ds.$
\end{thm}

The study of anisotropic operator is quite active in recent years.
One can find some of the interesting results about this topic, e. g.
in \cite{CGS, CS, CFV, DG1, DG2, DG3, FV, WX}, etc.

It is known that for any bounded smooth domain $\om$ in a complete
Riemannian manifold, there exists a unique solution $u_{\om}$,
called the torsion function of $\om$, to the equation
\be\label{th3eq1} \D u = -1 \ \ {\rm in} \ \ \om, \ \
u|_{\pa\om}=0\en and (\ref{th3eq1}) is the Euler equation of the
minimum problem \be\label{min} \min_{u\in W_0^{1,2}(\om)}
\int_{\om}\left(\fr 12|\na u|^2 -u\right). \en The number \be
\label{torsion} T(\om)=\int_{\om} u_{\om}= \int_{\om} |\na
u_{\om}|^2 \en is called the {\it torsional rigidity} of $\om$.

Serrin's theorem above  says that if the torsion function of a
bounded smooth domain  $\om$ in a Euclidean space has constant
derivative in the direction of the outward unit normal of $\pa\om$,
then $\om$ is a ball. In the third part of this paper, we give some
sharp estimates for the torsion function of a compact manifold with
boundary.

\begin{thm}\label{th3}
Let   $M$ be an $n$-dimensional compact Riemannian manifold with
boundary. Denote by  $T(M)$, $\rho$ and $V$ the torsion, the torsion
function and the volume of $M$, respectively.  Let $\nu$ be the
outward unit normal of $\pa M$ and assume that the Ricci curvature
of $M$ is bounded below by $(n-1)\kappa$.

i) We have \be \label{th3eq2} \min_{x\in \pa M}\fr{\pa^2 \rho}{\pa
\nu^2}(x)\leq -\fr 1n -\fr{(n-1)\kappa T(M)}V, \en with equality
holding if and only if $M$ is isometric to a ball in $\mathbb{R}^n$,
$\kappa=0$ and \be\label{th1eq00}\fr{\pa^2 \rho}{\pa \nu^2}=-\fr 1n\
\ {\rm on \ \ } \pa M. \en

ii) Let $A$ and $H$ be the area and the mean curvature of $\pa M$,
respectively. If $H\geq 0$ on $\pa M$,  then
 \be \label{th3eq3} \int_{\pa M} \fr{\pa^2
\rho}{\pa \nu^2}\leq (n-1)\left(\fr Vn -\kappa T(M)\right)^{\fr
12}\left(\int_{\pa M} H\right)^{\fr 12} -A, \en with equality
holding if and only if $\kappa=0$ and  $M$ is isometric to a ball
in $\mathbb{R}^n$.
\end{thm}
\begin{thm}\label{th4}
Let   $M$ be an $n$-dimensional compact Riemannian manifold with
boundary and  Ricci curvature bounded below by $(n-1)\kappa $.
 Let $u$ be the solution of the Dirichlet problem
\be\label{th4eq0}
  \D u=-1 \ \ {\rm in}\ \ \  M, \ \ u|_{\pa M}=0.
\en Then \be\label{th4eq1} \max_{\pa M}|\na u|^2\geq
\fr{(n+2)T(M)}{nV}+\fr{2(n-1)\kappa}V \int_{M} u|\na u|^2,\en with
equality holding if and only if $\kappa=0$ and $M$ is isometric to a
ball in $\mathbb{R}^n$.
\end{thm}

{\bf Remark 1.}  One can obtain Theorem B from Theorem \ref{th4}. In
fact, when $\om$ is a bounded smooth domain in $\mathbb{R}^n$, if
$u$ is a solution to the equation (\ref{serrin1}), then we have from
(\ref{inteq4}) that
$$c^2=\fr{(n+2)T(\om)}{nV}.$$
Thus, the equality sign in (\ref{th4eq1}) is attained  since the
Ricci curvature of $\om$ is zero. Theorem B follows.

The next result is a Pohozaev-type inequality on compact Riemannian
manifolds.

\begin{thm}\label{th6}
Let   $M$ be an $n$-dimensional compact Riemannian manifold with or
without  boundary. Assume that the  Ricci curvature of $M$ is
bounded below by $(n-1)\kappa $ and let $g: \mathbb{R}\ri\mathbb{R}$
be a continuous function.
 If  $u\in C^3(M)\cap C^1(\pa M)$ is a non-negative solution of the problem
\be\label{th6eq0}
  -\D u= g(u) \ \ {\rm in}\ \ \  M, \ \ u|_{\pa M}=0,
\en then we have \be\label{th6eq2} & & \int_M
g(u)\left(\fr{2(n-1)ug(u)}n -3G(u)-(n-1)\kappa u^2\right)\\ \no & &
\geq \left\{\begin{array}{l} \int_{\pa M} \left(\fr{\pa u}{\pa
\nu}\right)^3,\ \
{\rm when\ } \pa M\neq \emptyset, \\
0,  \ \ {\rm when\ } \pa M =\emptyset\end{array}\right.,\en where
$G(u)=\int_0^u g(\sigma)d\sigma.$
\end{thm}
From Theorem \ref{th6}, we have the following non-existence result.
\begin{cor} Let   $M$ be an $n$-dimensional closed Riemannian manifold
with Ricci curvature bounded below by $(n-1)\kappa $.  Assume that
$g: \mathbb{R}\ri\mathbb{R}$ is a continuous function and there
exists a discrete subset $S$ of $[0, +\infty)$ such that \be\no
g(t)\left(\fr{2(n-1)tg(t)}n -3\int_0^t g(\sigma)d\sigma-(n-1)\kappa
t^2\right)\left\{\begin{array}{l} =0, \ \ \ {\rm if } \ \ t\in S,
\\
<0, \ \ \ {\rm if } \ \ t\in [0, +\infty)\setminus S,
\end{array}\right..\en Then any non-negative solution of the equation \be \D u=- g(u)
\ \ {\rm on} \ \ M.\en is a constant.
\end{cor}

\section{Proof of Theorem \ref{th1} and Corollary \ref{cor2} }
In this section, we shall prove Theorem \ref{th1} and Corollary
\ref{cor2}. \vskip0.2cm {\it Proof of Theorem \ref{th1}.}
Multiplying the equation \be - \fr 1p {\rm
div}\left(|x|^{-bp}\na_{\xi}[F^p](\na u)\right)  = g(x, u) \en by
$-p\langle x, \na u\rangle$ and integrating on $\om$, one deduces
from divergence theorem that \be\label{th2peq1} & & -p\int_{\om}
g(x, u)\langle x, \na u\rangle dx
\\ \no &=& \int_{\om}{\rm
div}\left(|x|^{-bp}\na_{\xi}[F^p](\na u)\right)\langle x, \na
u\rangle dx \\ \no &=& \int_{\om}\left({\rm
div}\left(|x|^{-bp}\na_{\xi}[F^p](\na u)\langle x, \na
u\rangle\right)\right.\\ \no & & \ \ \ \ \ \ \ \ \ \  \
-\left.\left\langle |x|^{-bp}\na_{\xi}[F^p](\na u), \na\langle x,
\na u\rangle\right\rangle\right)dx
\\ \no &=& \int_{\pa\om}\left\langle |x|^{-bp}\na_{\xi}[F^p](\na u), \nu\right\rangle\langle x, \na u\rangle
d \mathcal{H}^{n-1}\\ \no & & -\int_{\om}\left\langle
|x|^{-ap}\na_{\xi}[F^p](\na u), \na u + \na^2 u(x)\right\rangle dx.
\en Here $\na^2 u : \mathfrak{X}(\om)\ri \mathfrak{X}(\om)$ denotes
the self-adjoint linear operator metrically equivalent to the
Hessian of $u$, and  is given by \cite{doC} \be\label{th2peq2}
\langle\na^2u(Z), W\rangle =\na^2u(Z, W)=\langle \na_Z\na u,
W\rangle\en for all $Z, W\in \mathfrak{X}(\om).$ It follows from
(\ref{aniseq3}) and
 $u|_{\pa\om}=0$ that
\be\label{th2peq3} & & \int_{\pa\om}\left\langle
|x|^{-bp}\na_{\xi}[F^p](\na u), \nu\right\rangle\langle x, \na
u\rangle d \mathcal{H}^{n-1}\\ \no &=& \int_{\pa\om}\left\langle
|x|^{-bp}\na_{\xi}[F^p](\na u),
\nu\right\rangle\langle x, u_{\nu} \nu\rangle d \mathcal{H}^{n-1}\\
\no &=& \int_{\pa\om}\left\langle |x|^{-bp}\na_{\xi}[F^p](\na u),
u_{\nu}\nu\right\rangle\langle x, \nu\rangle d \mathcal{H}^{n-1}
\\
\no &=& \int_{\pa\om}\left\langle |x|^{-bp}\na_{\xi}[F^p](\na u),
\na u\right\rangle\langle x, \nu\rangle d \mathcal{H}^{n-1}.
 \\
\no &=& p\int_{\pa\om}|x|^{-bp}F^p(\na u)\langle x, \nu\rangle d
\mathcal{H}^{n-1}.\en Using divergence theorem again, we infer
 \be
\label{th2peq4} -\int_{\om}\left\langle |x|^{-bp}\na_{\xi}[F^p](\na
u), \na u\right\rangle dx&=& \int_{\om} u\ {\rm div}(
|x|^{-bp}\na_{\xi}[F^p](\na u))dx\\
\no &=& -p\int_{\om} u g(x, u) dx. \en Let $\{e_1=(1,0,...,0),...,
e_n=(0,...,0, 1)\}$ be the canonical base of $\mathbb{R}^n$ and set
$u_i=\fr{\pa u}{\pa x_i}, \ u_{ij}=\fr{\pa^2 u}{\pa x_i x_j}, i,
j=1,...,n$.  We calculate   \be \left\langle\na_{\xi}[F^p](\na u),
\na^2 u(x)\right\rangle  &=& \left\langle\na_{\na_{\xi}[F^p](\na
u)}\na u, x\right\rangle
\\ \no &=& \left\langle\na_{\sum_{i=1}^n \fr{\pa [F^p]}{\pa \xi_i}(\na u)e_i}\na u, x\right\rangle
\\ \no &=& \sum_{i=1}^n \fr{\pa [F^p]}{\pa \xi_i}(\na u)\left\langle\na_{e_i}\na u, x\right\rangle
\\ \no &=& \sum_{i=1}^n \fr{\pa [F^p]}{\pa \xi_i}(\na u)\na^2 u(e_i,  x)
\\ \no &=& \sum_{i=1}^n \fr{\pa [F^p]}{\pa \xi_i}(\na u)\sum_{j=1}^n \langle x, e_j\rangle\na^2 u(e_i,  e_j)
\\ \no &=& \sum_{i, j=1}^n\fr{\pa [F^p]}{\pa \xi_i}(\na
u)\langle x, e_j\rangle  u_{ij}
\\ \no &=& \left\langle x, \sum_{i, j=1}^n\fr{\pa [F^p]}{\pa \xi_i}(\na
u)u_{ij} e_j\right\rangle
\\ \no &=& \left\langle x, \sum_{j=1}^n\fr{\pa (F^p(\na
u))}{\pa x_j}e_j\right\rangle \\ \no &=& \langle x, \na (F^p(\na
u))\rangle.
 \en
Therefore, we have \be \label{th2peq5}& &
\left\langle|x|^{-bp}\na_{\xi}[F^p](\na u), \na^2
u(x)\right\rangle\\ \no
 &=& |x|^{-bp}\langle\na(F^p(\na u)),
x\rangle\\ \no &=& \langle\na(F^p(\na u)), |x|^{-bp} x\rangle
\\ \no
&=& {\rm div}(F^p(\na u)|x|^{-bp}x) -F^p(\na u){\rm div}(|x|^{-bp}x)
\\ \no
&=& {\rm div}(F^p(\na u)|x|^{-bp}x) -(n-bp)|x|^{-ap}F^p(\na u).\en
Thus \be\label{th2peq6}& & \int_{\om}\left\langle
|x|^{-bp}\na_{\xi}F^p(\na u), \na^2 u(x)\right\rangle dx
\\
\no &=& \int_{\pa\om} |x|^{-bp}F^p(\na u)\langle x, \nu\rangle d
\mathcal{H}^{n-1} -(n-bp) \int_{\om} |x|^{-bp} F^p(\na u) dx
\\
\no &=& \int_{\pa\om} |x|^{-bp}F^p(\na u)\langle x, \nu\rangle d
\mathcal{H}^{n-1} -(n-bp) \int_{\om} |x|^{-bp}\cdot\fr 1p
\langle\na_{\xi}[F^p](\na u), \na u\rangle dx
\\ \no &=& \int_{\pa\om} |x|^{-bp}F^p(\na u)\langle x, \nu\rangle d \mathcal{H}^{n-1}+\fr{n-bp}p\int_{\om}u\ {\rm
div}\left(|x|^{-bp}\na_{\xi}[F^p](\na u)\right)
 dx
 \\ \no &=& \int_{\pa\om} |x|^{-bp}F^p(\na u)\langle x, \nu\rangle d \mathcal{H}^{n-1} -(n-bp)\int_{\om}u g(x, u) dx.\en
Substituting (\ref{th2peq3}), (\ref{th2peq4}) and (\ref{th2peq6})
 into (\ref{th2peq1}), we get
\be& &\label{th2peq7} -\int_{\om} g(x,u)\langle x, \na u\rangle dx\\
\no &=& \left(1- \fr 1p\right)\int_{\pa\om}|x|^{-bp} F^p(\na
u)\langle x, \nu\rangle d \mathcal{H}^{n-1}+\left(\fr np
-1-b\right)\int_{\om} u g(x, u) dx.
 \en
On the other hand, we have \be\label{th2peq8} g(x, u)\langle x, \na
u\rangle =\langle x, \na(G(x, u))\rangle - \langle x, \na_x G(x,
u)\rangle \en and so \be \label{th2peq9}  -\int_{\om} g(x, u)\langle
x, \na u\rangle dx &=& -\int_{\om} \left(\langle x, \na(G(x,
u))\rangle - \langle x, \na_x G(x,
u)\rangle\right) dx\\
\no &=& n\int_{\om} G(x, u) dx +\int_{\om}\langle x, \na_xG(x,
u)\rangle dx. \en Combining (\ref{th2peq7}) and
(\ref{th2peq9}), we get (\ref{th1eq2}). \qed\\

\vskip0.3cm {\it Proof of Corollary 1.3.} Suppose that $u$ is a
positive solution of (\ref{th1eq4}). Then the equality
(\ref{th1eq2}) holds. Thus, we have
 \be \label{cor1eq4} & &\left(1+b-\fr np\right)\int_{\om} u (\la |x|^{-\alpha}u^{r-1}+\mu |x|^{-\beta} u^{s-1}
+\eta |x|^{-\gamma}u^{t-1}) dx\\ \no & & \ \ +n\int_{\om} G(x, u) dx
+\int_{\om}\langle x, \na_xG(x, u)\rangle dx\\ \no &=& \left(1- \fr
1p\right)\int_{\pa\om} F^p(\na u)\langle x, \nu\rangle ds\\ \no
&\geq & 0, \en where \be \label{cor1eq5} G(x, u)&=&\int_0^u (\la
|x|^{-\alpha}|\sigma|^{r-2}\sigma+\mu |x|^{-\beta}
|\sigma|^{s-2}\sigma
+\eta |x|^{-\gamma}|\sigma|^{t-2}\sigma)d\sigma\\
\no &=&\fr{\la}r |x|^{-\alpha} u^r+ \fr{\mu}s
|x|^{-\beta} u^s +\fr{\eta}t |x|^{-\beta} u^t, \en \be\label{cor1eq6} & & \langle x, \na_x G(x, u)\\
\no &=& \sum_{i=1}^n x_i\fr{\pa G}{\pa x_i}(x, u)\\ \no &=&
-\fr{\la\alpha}r \alpha |x|^{-\alpha}u^r-\fr{\mu\beta}s |x|^{-\beta}
u^s-\fr{\eta\gamma}t |x|^{-\gamma} u^t. \en Substituting
(\ref{cor1eq5}) and (\ref{cor1eq6}) into (\ref{cor1eq4}), we have
\be& &
 \la\left(1+b-\fr np
+\fr{n-\alpha}r\right)\int_{\om}|x|^{-\alpha}u^r dx+\mu\left(
1+b-\fr np+\fr{n-\beta}s\right)\int_{\om}  |x|^{-\beta} u^s dx \\
\no & &  + \eta\left(1+b-\fr np+\fr{n-\gamma}t\right)\int_{\om}
|x|^{-\gamma} u^t dx\geq 0. \en This is a contradiction if
(\ref{cor2eq1}) holds.
\qed\\
Using Theorem \ref{th1} we can also prove the following nonexistence
result.
\begin{cor}\label{cor3}
Let $\om$ be a bounded star-shaped domain with smooth boundary and
$b, \alpha, \beta, \lambda, \mu,  r$ be constants such that $ r\neq
0$ and \be\label{cor3eq0} \la\left(n+1+b-\fr np-\alpha\right)\leq 0,
\ \ \mu\left( 1+b-\fr np+\fr{n-\beta}r\right)<0. \en Then, the
problem \be\label{cor3eq1}\left\{\begin{array}{l} - \fr 1p {\rm
div}\left(|x|^{-bp}\na_{\xi}[F^p](\na u)\right)  = \la
|x|^{-\alpha}+\mu |x|^{-\beta} |u|^{r-2} u \  \ \ {\rm in\ \
 }\om, \\
u|_{\pa\om}=0,
\end{array}\right.
\en has no positive solution.
\end{cor}
{\it Proof of Corrolary \ref{cor3}.} If $u$ is a positive solution
of (\ref{cor3eq1}), then we have from (\ref{th1eq2}) that
 \be \label{cor3eq4} & &\left(1+b-\fr np\right)\int_{\om} u (\la |x|^{-\alpha}+\mu |x|^{-\beta} u^{r-1}) dx
 \\ \no & & \ \ +n\int_{\om} G(x, u) dx +\int_{\om}\langle x,
\na_xG(x, u)\rangle dx\\ \no &=& \left(1- \fr 1p\right)\int_{\pa\om}
F^p(\na u)\langle x, \nu\rangle ds\\ \no &\geq & 0. \en Here \be
\label{cor3eq5} G(x, u)&=&\int_0^u (\la |x|^{-\alpha}+\mu
|x|^{-\beta} |\sigma|^{r-2}\sigma)d\sigma\\
\no &=&\la |x|^{-\alpha} u + \fr{\mu}r |x|^{-\beta} u^r, \en
\be\label{cor3eq6} & &\langle x, \na_x G(x, u)\\
\no &=& \sum_{i=1}^n x_i\fr{\pa G}{\pa x_i}(x, u)\\ \no &=&
-\la\alpha  |x|^{-\alpha}u-\fr{\mu\beta}r |x|^{-\beta} u^r. \en
Substituting (\ref{cor3eq5}) and (\ref{cor3eq6}) into
(\ref{cor3eq4}), we have \be\no
 \la\left(n+1+b-\fr np-\alpha\right)\int_{\om}|x|^{-\alpha}u dx+\mu\left(
1+b-\fr np+\fr{n-\beta}r\right)\int_{\om}  |x|^{-\beta} u^r dx \geq
0, \en contradicting to (\ref{cor3eq0}).
\qed\\

\section{Proof of Theorems \ref{th2} and \ref{th0}}
In this section we prove Theorems \ref{th2} and \ref{th0}. Firstly
we recall some facts needed about the function $F$ introduced in
section 1. Because of (\ref{ineq1}) we can assume, without loss of
generality, that the convex closed set
$$ K=\{ x\in \mathbb{R}^n: F(x)\leq 1\}$$
has measure $|K|$ equal to the measure $\omega_n$ of the unit sphere
in $\mathbb{R}^n$. We say that $F$ is the gauge of $K$. The support
function of $K$ is defined as \cite{R} \be F^o(x)=\sup_{\xi\in
K}\langle x, \xi\rangle. \en It is easy to see that $F^o :
\mathbb{R}^n\ri [0, +\infty)$ is a convex, homogeneous function and
that $F, F^o$ are polar each other in the sense that \be
F^o(x)=\sup_{\xi\neq 0}\fr{\langle x, \xi\rangle}{F(\xi)}, \en and
\be F(x)=\sup_{\xi\neq 0}\fr{\langle x, \xi\rangle}{F^o(\xi)}. \en
We set
$$ K^o=\{ x\in \mathbb{R}^n: F^o(x)\leq 1\}$$
and denote by $\kappa_n$ the measure of $K^o$.

Let $\om$ be an open subset of $\mathbb{R}^n$. The total variation
of a function $u\in BV(\om)$ with respect to a guage fuction $F$ is
defined by \cite{AB} \be \int_{\om} |\na u|_{F}
=\sup\left\{\int_{\om} u\ {\rm div}\ \sigma\ dx: \sigma\in
C_0^1(\om; \mathbb{R}^n), F^o(\sigma)\leq 1\right\}. \en The
perimeter of a set $E$ with respect to $F$ is then defined as \be
P_F(E; \om)=\int_{\om} |\na \chi_E|_{F} =\sup\left\{\int_{\om}  div\
\sigma\ dx: \sigma\in C_0^1(\om; \mathbb{R}^n), F^o(\sigma)\leq
1\right\}. \en The following co-area formula \be \int_{\om}|\na
u|_F=\int_0^{\infty}P_F(\{u>s\}; \om) ds, \ \ \forall u\in BV(\om),
\en and the equality \be P_F(E; \om)=\int_{\om\cap \pa^*E} F(\nu^E)
d{\mathcal H}^{n-1} \en hold, where $\pa^*E$ is the reduced boundary
of $E$ and $\nu^E$ is the outer normal to $E$ (see \cite{AB}). The
following result can be found in \cite{AFTL}, \cite{DP}, \cite{FM}.
\begin{lem}\label{lem1}(Wulff theorem). If $E$ is a set of finite
perimeter in $\mathbb{R}^n$, then \be\label{lem2} P_F(E;
\mathbb{R}^n)\geq n\kappa_n^{1/n}|E|^{1-1/n}, \en and equality holds
if and only if $E$ has Wulff shape, i.e., E is a sub-level set of
$F^o$, modulo translations.
\end{lem}
Now we are ready to give a

 {\it Proof of Theorem \ref{th2}}. Let
$v_{\om}$ be the unique solution of the equation \be\label{pth2eq1}
-\mathcal{Q}_p v=1, \ \ \ {\rm in}\ \ \ \om, \ \ \ \ v=0 \ \ \ {\rm
on }\ \ \ \pa \om. \en Then $v_{\om}$ is positive in $\om$. By
(\ref{th201}) and since $F\in C^3(\mathbb{R}^n\setminus\{0\})$, we
know that $v_{\om}\in C^{1, \alpha}(\om) \cap C^3((\{\na v_{\om}\neq
0\})$ (see \cite{LU, To}).

It follows from (\ref{th2eq1}) that \be\label{pth2eq2}
\la_{p,1}(F,\om)\leq \fr{\int_{\om} (F(\na v_{\om}))^p
dx}{\int_{\om} v_{\om}^p dx}= \fr{\int_{\om} v_{\om} dx}{\int_{\om}
v_{\om}^p dx}. \en Combining (\ref{th2eq4}) and (\ref{pth2eq2}), we
infer \be\label{pth2eq3} \la_{p,1}(F,\om)T_{F,p}(\om)^{p-1}\leq
\fr{\left(\int_{\om} v_{\om} dx\right)^p}{\int_{\om} v_{\om}^p dx}.
\en
 Let $M=\sup_{\om}
v_{\om}$. For $s\in [0, M]$, we denote by \be\label{pth2eq4}
\mu(s)=|\{x\in \om : v_{\om}> s\}| \en the distribution function of
$v_{\om}$. Then \be\label{pth2eq4} \int_{\om} v_{\om} =\int_0^M
\mu(s) ds \en and \be\label{pth2eq5} \int_{\om} v_{\om}^p dx
=\int_0^M ps^{p-1} \mu(s)ds. \en Observe that the boundary of
\be\label{pth2eq6} \{x\in \om : v_{\om}> s\} \en
 is
\be\label{pth2eq7} \{x\in \om : v_{\om}= s\} \en for almost every
$s>0$ and the inner normal to this boundary at a point $x$ is
exactly $\na v_{\om}(x)/|\na v_{\om}(x)|$. Integrating
$-\mathcal{Q}_p v_{\om}=1$ over (\ref{pth2eq6}) gives \be
\label{pth2eq8} \mu(s)&=& -\fr 1p\int_{v_{\om}(x)>s} {\rm
div}\left(\na_{\xi}[F^p](\na v_{\om})\right) dx\\ \no &=& \fr
1p\int_{v_{\om}(x)=s}\left\langle \na_{\xi}[F^p](\na v_{\om}),
\fr{\na v_{\om}}{|\na v_{\om}|}\right\rangle d\mathcal{H}^{n-1}
\\ \no &=&
\int_{v_{\om}(x)=s}\fr{F^p(\na v_{\om})}{|\na v_{\om}|}
d\mathcal{H}^{n-1} \en and we have
 \be \label{pth2eq9}
-\mu^{\prime}(s)=\int_{v_{\om}(x)=s}\fr 1{ |\na v_{\om}|} d{\mathcal
H}^{n-1} \en for almost every $s\in [0, M)$.

The co-area formula gives that \be \label{pth2eq10}
-\fr{d}{dt}\int_{v_{\om}>s} F(\na v_{\om}) dx=P_{F}(v_{\om}>s\};
\om), \en for almost all $s$. Also, since $v_{\om}$ is smooth with
compact support, it is known \cite{AFTL} that for almost  every
$s\in [0, M),$
 \be \label{pth2eq10}
-\fr{d}{dt}\int_{u>s} F(\na v_{\om}) dx=\int_{v_{\om}=t}\fr{F(\na
v_{\om})}{|\na v_{\om}|} d\mathcal{H}^{n-1}. \en Hence, \be
\label{pth2eq11} P_{F}(\{v_{\om}>s\}; \om)=\int_{v_{\om}=s}\fr{F(\na
v_{\om})}{|\na v_{\om}|} d\mathcal{H}^{n-1}. \en  From H\"older's
inequality, (\ref{pth2eq8}) and (\ref{pth2eq9}), we obtain
\be\label{pth2eq12}\left(P_{F}(\{v_{\om}>s\}; \om)\right)^p&=&
\left(\int_{v_{\om}=s}\fr{F(\na v_{\om})}{|\na v_{\om}|} d\mathcal{H}^{n-1}\right)^p\\
\no &\leq& \int_{v_{\om}=s}\fr{F(\na v_{\om})^p}{|\na v_{\om}|}
d\mathcal{H}^{n-1}\left(\int_{v_{\om}=s}\fr 1{|\na
v_{\om}|}d\mathcal{H}^{n-1}\right)^{p-1} \\ \no &=&
\mu(s)(-\mu^{\prime}(s))^{p-1}. \en The isoperimetric inequality
(\ref{lem2}) tells us that \be\label{pth2eq13} P(\{v_{\om}>s\})\geq
n\kappa_n^{1/n}\mu(s)^{(n-1)/n}. \en One can then use the same
arguments as in the proof of Theorem 1 in \cite{xia} to finish the
proof of Theorem \ref{th2}. For the sake of completeness, we include
it. It follows from (\ref{pth2eq12}) and (\ref{pth2eq13}) that
\be\label{pth2eq14} \left(n\kappa_n^{1/n}\right)^{\fr p{p-1}}\leq
-\mu(s)^{\fr{n+p-np}{n(p-1)}} \mu^{\prime}(s) \en Integrating
(\ref{pth2eq14}) gives \be\label{pth2eq15} \mu(s)&\leq&
\left(|\om|^{\fr{p}{n(p-1)}}-\left(n\kappa_n^{1/n}\right)^{\fr
p{p-1}}\cdot\fr{p}{n(p-1)}s\right)^{\fr{n(p-1)}p}\\
\no &= & |\om|(1-bs)^a, \en where \be\label{pth2eq16}
a=\fr{n(p-1)}{p}, \ \ b=\left(n\kappa_n^{1/n}\right)^{\fr
p{p-1}}\cdot\fr{p}{n(p-1)}\cdot|\om|^{-\fr 1a}. \en Define $F: [0,
M]\rightarrow {\mathbb R}$ by \be\label{pth2eq17}
F(t)=\left(\int_0^t \mu(s) ds\right)^p - p\left(\int_0^t
s^{p-1}\mu(s) ds\right)|\om|^{p-1}. \en It is easy to see from
(\ref{pth2eq15}) that \be\label{pth2eq18} F^{\prime}(t) &=&
p\left(\left(\int_0^t \mu(s) ds\right)^{p-1} - t^{p-1}
|\om|^{p-1}\right)\mu(t)\\ \no &\leq& p|\om|^{p-1}\left(\left(\fr
1{b(a+1)}\left(1-(1-bt)^{a+1}\right)\right)^{p-1}-t^{p-1}\right)\mu(t)
\en Since $p\geq 2\geq n/(n-1)$, we have $a+1\geq 2$. Using  \be
(1+x)^{\alpha}\geq 1+\alpha x +x^2, \ \alpha \geq 2,\ x\geq -1 \en
in (\ref{pth2eq18}), we obtain
 \be\label{pth2eq19}
F^{\prime}(t)&\leq&p|\om|^{p-1}
\left(\left(t-\fr {bt^2}{a+1}\right)^{p-1}-t^{p-1}\right)\mu(t)\\
\no &=&p|\om|^{p-1}t^{p-1} \left(\left(1-\fr
{bt}{a+1}\right)^{p-1}-1\right)\mu(t)\\ \no &\leq &
-\fr{pb|\om|^{p-1}t^p\mu(t)} {a+1} \en Integrating (\ref{pth2eq19})
over $[0, M]$ and using H\"older's inequality we have \be F(M)&\leq&
-\fr{pb|\om|^{p-1} } {a+1} \int_0^M t^p\mu(t)dt
\\ \no &\leq& -\fr{pb|\om|^{p-1}}{a+1}\cdot\fr{\left(\int_0^M
t^{p-1}\mu(t)dt\right)^{p/(p-1)}}{\left(\int_0^M
\mu(t)dt\right)^{1/(p-1)}}\\ \no &=& -\fr{b|\om|^{p-1} } {a+1}
\cdot\fr{\left(\int_{\om} v_{\om}^pdx\right)^{p/(p-1)}}{p^{\fr
1{p-1}}\left(\int_{\om} v_{\om}dx\right)^{1/(p-1)}}, \en that is,
\be\label{pth2eq20} \left(\int_{\om} v_{\om}dx\right)^p -
|\om|^{p-1}\int_{\om} v_{\om}^p dx \leq -\fr{b|\om|^{p-1} } {a+1}
\cdot\fr{\left(\int_{\om} v_{\om}^pdx\right)^{p/(p-1)}}{p^{\fr
1{p-1}}\left(\int_{\om} v_{\om}dx\right)^{1/(p-1)}}. \en Dividing by
$|\om|^{p-1}\int_{\om} v_{\om}^p dx$ and using H\"older's
inequality, we infer
 \be\label{pth2eq21}
 & &\fr{\la_{p,1}(F,\om)T_{F,p}(\om)^{p-1}}{|\om|^{p-1}}-1\\ \no &\leq &
\fr{\left(\int_{\om} v_{\om}dx\right)^p}{\left(\int_{\om}
v_{\om}^pdx\right) |\om|^{p-1}} -1\\ \no &\leq& -\fr b{(a+1)p^{\fr
1{p-1}}}\cdot\left(\fr{\int_{\om} v_{\om}^pdx}{\int_{\om}
v_{\om}dx}\right)^{1/(p-1)}\\ \no &\leq&
 -\fr b{(a+1)p^{\fr 1{p-1}}}
 \cdot\fr{\int_{\om} v_{\om}dx}{|\om|}\\ \no &=&
- \fr{p^{\fr{2p-3}{p-1}}(n\kappa_n^{1/n})^{\fr
p{p-1}}}{n(p-1)(n(p-1)+p)}\cdot \fr{T_{F,p}(\om)}{|\om|^{1+\fr
p{n(p-1)}}}.
 \en
Thus (\ref{th2eq7}) holds. \qed\\

{\it Proof of Theorem \ref{th0}.} Let $u_M=\sup_{\om} u$ and
consider two functions $\eta, V: [0, u_M]\ri \mathbb{R}$ given by
\be\label{pth0eq1} \eta(t)=\int_{\{x\in \om: u(x)>t\}} g(u) dx, \ \
V(t)=|\{x\in \om: u(x)>t\}|.\en Integrating the equation $-\fr 1n
{\rm div}(\na_{\xi}[F^n](\na u))= g(u)$ on $\{x\in \om: u(x)>t\}$,
we have \be\label{pth0eq2} \eta(t) = \int_{\Gamma(t)} \fr{F^n(\na
u)}{|\na u|} d\mathcal{H}^{n-1}, \en where $\Gamma(t)=\{ x\in \om :
u(x)=t\}.$ The co-area formula gives
 \be\label{pth0eq3}
-\fr{dV}{dt}=\int_{\{x\in \om: u(x)=t\}}\fr{d\mathcal{H}^{n-1}}{|\na
u|}. \en From (\ref{pth0eq2}), (\ref{pth0eq3}) and H\"older and
inequalities we get \be\label{pth0eq4} \eta (-V^{\prime})^{n-1}&=&
\int_{\Gamma(t)}\fr{F^n(\na u)}{|\na u|}
d\mathcal{H}^{n-1}\left(\int_{\{x\in \om:
u(x)=t\}}\fr{d\mathcal{H}^{n-1}}{|\na u|}\right)^{n-1}\\ \no &\geq&
\left(\int_{\Gamma(t)}\fr{F(\na u)}{|\na u|}
d\mathcal{H}^{n-1}\right)^n\\ \no &=& \left(P_F(\{ u>t\}; \om)\right)^n\\
\no &\geq& \left(n\kappa_n^{\fr 1n}V(t)^{\fr{n-1}n}\right)^n .\en
Hence \be -\eta^{\fr 1{n-1}}V^{\prime}(t)\geq
\left(n^n\kappa_n\right)^{\fr 1{n-1}} V(t),\en which, combining with
\be \fr{d\eta}{dt}=g(t) \fr{dV}{dt}, \en gives \be \label{pth0eq4}
-\eta^{\fr 1{n-1}}\eta^{\prime}(t)\geq \left(n^n\kappa_n\right)^{\fr
1{n-1}} g(t)V(t).\en Integration of (\ref{pth0eq4}) on $[0, u_M]$
yields \be\label{pth0eq5}\fr{n-1}n \eta(0)^{\fr n{n-1}}&\geq&
\left(n^n\kappa_n\right)^{\fr 1{n-1}}\int_0^{u_M} g(t)V(t)dt\\ \no
&=& -\left(n^n\kappa_n\right)^{\fr 1{n-1}}\int_0^{u_M}
G(t)V^{\prime}(t)dt\\ \no &=& \left(n^n\kappa_n\right)^{\fr
1{n-1}}\int_{\om} G(u) dx,\en
 that is
 \be \fr{n-1}n\left(\int_{\om} g(u) dx\right)^{\fr n{n-1}}\geq
\left(n^n\kappa_n\right)^{\fr 1{n-1}}\int_{\om} G(u) dx.
 \en
(\ref{th0e2}) follows. \qed\\

\section{Proof of Theorems \ref{th3}-\ref{th6} }

In this section, we prove Theorems
\ref{th3} and \ref{th4}.
 Before doing this, we first recall Reilly's formula which will be used later. Let $M$
be an $n$-dimensional compact manifold with boundary. We will often
write $\langle, \rangle$ the Riemannian metric on $M$ as well as
that induced on $\pa M$. Let $\na$ and $\D $ be the connection and
the Laplacian on $M$, respectively. Let $\nu$ be the unit outward
normal vector of $\pa M$. The shape operator of $\pa M$ is given by
$S(X)=\na_X \nu$ and the second fundamental form of $\pa M$ is
defined as $II(X, Y)=\langle S(X), Y\rangle$, here $X, Y\in T (\pa
M)$. The eigenvalues of $S$ are called the principal curvatures of
$\pa M$ and the mean curvature $H$ of $\pa M$ is given by $H=\fr
1{n-1} {\rm tr\ } S$, here ${\rm tr\ } S$ denotes the trace of $S$.
 For a smooth function
$f$ defined on $M$, the following identity holds \cite{RE} if
$h=\left.\fr{\pa }{\pa \nu}f\right|_{\pa M}$, $z=f|_{\pa M}$ and
${\rm Ric}$ denotes the Ricci tensor of $M$: \be & & \int_M
\left((\D f)^2-|\na^2 f|^2-{\rm Ric}(\na f, \na f)\right)\\ \no &=&
\int_{\pa M}\left( ((n-1)Hh+2\overline{\D}z)h + II(\overline{\na}z,
\overline{\na}z)\right). \en Here $\na^2 f$ is the Hessian of $f$;
$\ov{\D}$ and $\ov{\na}$  represent the Laplacian and the gradient
on $\pa M$ with respect to the induced metric on $\pa M$,
respectively.
 \vskip0.2cm

{\it Proof of Theorem \ref{th3}.} $i)$ Since $\rho$ satisfies the
equation \be \label{pth3eq1} \D \rho = -1 \ \ {\rm in}\ \ M, \ \ \
\rho|_{\pa M}=0, \en we know from the strong maximum principle and
Hopf lemma \cite{GT} that $\rho$ is positive in the interior of $M$
and \be \label{pth3eq2}\fr{\pa \rho}{\pa \nu}(x)<0, \ \  \ \forall
x\in \pa M. \en It follows from Bochner formula that
\be\label{pth3eq3} \fr 12 \D |\na \rho|^2& =& |\na^2 \rho|^2
+\langle \na \rho, \na(\D \rho)\rangle +{\rm Ric}(\na \rho, \na
\rho)\\ \no & =& |\na^2 \rho|^2 + {\rm Ric}(\na \rho, \na \rho)
\\  \no &\geq & |\na^2 \rho|^2 +(n-1)\kappa|\na \rho|^2. \en
 Integrating
(\ref{pth3eq3}) on $M$ and using divergence theorem, we get
\be\label{pth3eq4}\int_M |\na^2 \rho|^2  +(n-1)\kappa\cdot
T(M)&=&\int_M ( |\na^2 \rho|^2+(n-1)\kappa|\na \rho|^2)\\ \no &\leq&
\fr 12\int_{\pa M} \nu |\na \rho|^2\\ \no &=& \int_{\pa M} \na^2
\rho(\na \rho, \nu). \en Since $\rho|_{\pa M}=0$, we have \be\no \na
\rho|_{\pa M}=\left(\fr{\pa \rho}{\pa \nu}\right)\nu, \ \ \na^2
\rho(\nu, \nu)= \fr{\pa^2 \rho}{\pa\nu^2}.\en Hence \be\no \int_{\pa
M} \na^2 \rho(\na\rho, \nu) = \int_{\pa M} \left(\fr{\pa \rho}{\pa
\nu}\right)\left(\fr{\pa^2 \rho}{\pa\nu^2}\right). \en Setting
$$ l= \min_{x\in \pa M}\fr{\pa^2 \rho}{\pa \nu^2}(x);$$
we have from (\ref{pth3eq2}) that \be\label{pth3eq5} \left(\fr{\pa
\rho}{\pa \nu}\right)\left(\fr{\pa^2 \rho}{\pa\nu^2}\right)\leq
\left(\fr{\pa \rho}{\pa \nu}\right)l.\en Hence \be\label{pth3eq6}
\int_{\pa M} \na^2 \rho(\na\rho, \nu) &\leq & l\int_{\pa M}\fr{\pa
\rho}{\pa \nu}
\\ \no &=& l\int_{ M} \D \rho\\ \no &=& -l V.
\en The Schwarz inequality implies that \be \label{pth3eq7} |\na^2
\rho|^2\geq \fr 1n (\D \rho)^2 =\fr 1n \en with equality holding if
and only if \be\label{pth3eq8} \na^2 \rho =\fr{\D \rho}n \langle ,
\rangle = - \fr 1n \langle , \rangle. \en Combining (\ref{pth3eq4}),
(\ref{pth3eq6}) and (\ref{pth3eq7}), we get (\ref{th3eq2}). On the
other hand, if (\ref{th3eq2}) take equality sign, then the
inequalities (\ref{pth3eq3})-(\ref{pth3eq7}) should be  equalities.
Thus,
 (\ref{pth3eq8}) holds on $M$. Taking the
covariant derivative of (\ref{pth3eq8}), we get
 $\na^3 \rho=0$ and from the Ricci identity,
\be \label{pth3eq9} R(X, Y)\na \rho=0, \en for any tangent vectors
$X, Y$ on $M$, where $R$ is the curvature tensor of $M$. By the the
maximum principle $\rho $ attains its maximum at some point $x_0$ in
the interior of $M$. Let $r$ be the distance function to $x_0$; then
from (\ref{pth3eq8}) it follows that \be\label{pth3eq10} \na \rho =
-\fr 1n r \fr{\pa}{\pa r}. \en Using (\ref{pth3eq9}),
(\ref{pth3eq10}), Cartan's theorem (cf. \cite{doC}) and $\rho|_{\pa
M}=0$, we conclude that $M$ is a ball in $\mathbb{R}^n$ whose center
is $x_0$, and \be \no \rho(x)=\fr 1{2n}(r_0^2 - |x-x_0|)\en in $M$,
here $r_0$ is the radius of the ball. This in turn implies that
$\kappa =0$.

$ii)$ Restricting $\D \rho = -1$ on $\pa M$ and noticing $\rho|_{\pa
M}=0$, we infer \be\label{pth3eq11}
\fr{\pa^2\rho}{\pa\nu^2}+(n-1)H\fr{\pa\rho}{\pa\nu} = -1 \ \ {\rm on
} \ \ \pa M. \en Integrating  (\ref{pth3eq11}) on $\pa M$ yields
\be\label{pth3eq12} A+ \int_{\pa M} \fr{\pa^2\rho}{\pa\nu^2} =
-(n-1)\int_{\pa M} H\pa_{\nu}\rho. \en Substituting $\rho$ into
Reilly formula, we get \be\label{pth3eq13} (n-1)\int_{\pa
M}H(\pa_{\nu}\rho)^2 &=&\int_M ((\D \rho)^2-|\na^2 \rho|^2-{\rm
Ric}(\na\rho, \na\rho)) \\ \no &\leq& \int_M \left((\D \rho)^2-\fr
1n (\D \rho)^2-(n-1)\kappa|\na\rho|^2\right) \\ \no &=&
\fr{(n-1)V}n-(n-1)\kappa T(M), \en with equality holding if and only
if \be \label{pth3eq14} |\na^2\rho|^2 =\fr 1n\ \ {\rm on\ \ }M,\en
and \be \label{pth3eq15} {\rm Ric}(\na\rho,
\na\rho)=(n-1)\kappa|\na\rho|^2 \ \ {\rm on\ \ } M. \en Since $H\geq
0$ on $\pa M$, one obtains from  H\"older's inequality that
\be\label{pth3eq16} -\int_{\pa M} H\pa_{\nu}\rho\leq \left(\int_{\pa
M} H(\pa_{\nu}\rho)^2\right)^{\fr 12}\left(\int_{\pa M}
H\right)^{\fr 12}. \en Combining (\ref{pth3eq12}), (\ref{pth3eq13})
and (\ref{pth3eq16}), one gets ((\ref{th3eq3}). Also, the equality
in (\ref{th3eq3}) holding implies that (\ref{pth3eq14}) holds. Using
the same arguments as in the proof of item $i)$, we conclude that
$M$ is isometric to a
ball in $\mathbb{R}^n$.\qed\\

{\it Proof of Theorem \ref{th4}.} As stated in the proof of Theorem
\ref{th3}, the function $u$ is positive in the interior of $M$ and
\be\label{th4eq2} -\left.\fr{\pa u}{\pa \nu}\right|_{\pa M} >0. \en
Multiplying the equation \be\label{th4eq3}
 -\D u= 1\en
  by $ \left(|\na u|^2 +\fr 2n
u\right)$ and integrating on $M$, we have from divergence theorem
and $u|_{\pa M}=0$ that \be\label{th4eq4}& & \fr{n+2}n T(M)\\ \no &=&\fr{n+2}n \int_{M}u \\
\no &=&\int_{M} \left( |\na u|^2 +\fr 2n u\right)\\ \no &=& \int_{M}
\left( |\na u|^2 +\fr 2n u\right)(-\D u)\\ \no &=&\int_M\left\langle
\na\left(|\na u|^2 +\fr 2n u\right), \na u\right\rangle -\int_{\pa
M}\left(|\na u|^2 +\fr 2n u\right)\fr{\pa u}{\pa\nu}\\ \no &=&
 -\int_{M} u\D \left( |\na u|^2 +\fr 2n
u\right)-\int_{\partial M} \left( |\na u|^2 +\fr 2n u\right)\fr{\pa
u}{\pa \nu}
\\ \no
&=& -2\int_{M} u\left( |\na^2 u|^2 +{\rm Ric}(\na u, \na u)-\fr 1n
\right)-\int_{\partial M} \left(\fr{\pa u}{\pa \nu}\right)^3
\\ \no
&\leq& -2\int_{M} u {\rm Ric}(\na u, \na u)-\int_{\partial M}
\left(\fr{\pa u}{\pa \nu}\right)^3,\en with equality holding if and
only if
 \be
\label{th4eq5} |\na^2 u|^2=\fr 1n \ \ {\rm on} \ \ M. \en Setting
\be\label{th4eq6} \max_{x\in\pa M}|\na u|=m\en and using \be
\label{th4eq6} u {\rm Ric}(\na u, \na u)\geq (n-1)\kappa u|\na
u|^2,\en we conclude from (\ref{th4eq4}) that \be \label{th4eq7}
\fr{n+2}n T(M)+2(n-1)\kappa \int_{M} u|\na
u|^2&\leq& m^2\int_{\partial M} \left(-\fr{\pa u}{\pa \nu}\right)\\
\no &=& m^2\int_M (-\D u)\\ \no &=& m^2 V. \en Thus (\ref{th4eq1})
holds. It is clear from the above proof that if the equality in
(\ref{th4eq1}) holds then (\ref{th4eq5}) holds and so $\kappa=0$ and
$M$ is isometric to a ball in
$\mathbb{R}^n$.  \qed\\

{\it Proof of Theorem \ref{th6}.} We shall only consider the case
that $\pa M\neq \emptyset$ since the case $\pa M=\emptyset$ is
similar. Multiplying the equation $\D u= - g(u)$ by $|\na u|^2$,
integrating on $M$ and using the divergence theorem, we get
\be\label{pth6eq1} & & \int_M g(u)|\na u|^2\\ \no &=&- \int_M |\na
u|^2 \D u\\ \no &=& \int_M \langle \na|\na u|^2, \na u\rangle
-\int_{\pa M} |\na u|^2 \fr{\pa u}{\pa \nu}\\ \no &=& -\int_M u\D
|\na u|^2 -\int_{\pa M} \left(\fr{\pa u}{\pa \nu}\right)^3\\ \no &=&
-\int_M 2u(|\na^2u|^2 +{\rm Ric}(\na u, \na u)+\langle \na u, \na(\D
u))\rangle) - \int_{\pa M} \left(\fr{\pa u}{\pa \nu}\right)^3\\ \no
&\leq& -\int_M 2u\left(\fr{(\D u)^2}n +(n-1)\kappa|\na u|^2+\langle
\na u, \na(\D u)\rangle\right) - \int_{\pa M} \left(\fr{\pa u}{\pa
\nu}\right)^3
\\ \no &=& -\int_M \left(\fr{2ug(u)^2}n
+(n-1)\kappa \langle \na u, \na u^2\rangle+\langle \na u^2, \na(\D
u)\rangle\right) - \int_{\pa M} \left(\fr{\pa u}{\pa \nu}\right)^3
\\ \no &=& \int_M \left(-\fr{2ug(u)^2}n
+(n-1)\kappa u^2\D u + \D u \D u^2\right) - \int_{\pa M}
\left(\fr{\pa u}{\pa \nu}\right)^3
\\ \no &=& \int_M \left(-\fr{2ug(u)^2}n
-(n-1)\kappa u^2g(u) + 2ug(u)^2-2g(u)|\na u|^2\right) - \int_{\pa M}
\left(\fr{\pa u}{\pa \nu}\right)^3.\en On the other hand, it is easy
to see that\be \label{pth6eq2} \int_M g(u)|\na u|^2 =\int_M \langle
\na G(u), \na u\rangle =-\int_M G(u) \D u = \int_M G(u) g(u). \en
Substituting (\ref{pth6eq2}) into (\ref{pth6eq1}), one gets
(\ref{th6eq2}).
\qed\\

\end{document}